\documentclass[11pt,english]{article}

\usepackage{anysize}
\marginsize{2cm}{2cm}{2.5cm}{2.4cm}

\usepackage[english]{babel}
\usepackage{amsfonts}
\usepackage{latexsym}
\usepackage{amsmath}
\usepackage{mathrsfs}
\usepackage{amssymb}
\usepackage{amsthm}
 \usepackage{graphicx}
\usepackage[usenames]{color}

\let\phi=\varphi
\newcommand{\Z}{{\mathbb Z}}

\newcommand{\N}{{\mathbb N}}
\newcommand{\eps}{\varepsilon}
\newcommand{\IP}{{\mathbb P}}
\newcommand{\IE}{{\mathbb E}}

\newcommand{\Nc}{\mathcal{N}}
\newcommand{\Ec}{\mathcal{E}}
\newcommand{\Ic}{\mathcal{T}}
\newcommand{\Ac}{\mathcal{A}}
\newcommand{\Wc}{\mathcal{W}}
\newcommand{\Par}{\o}
\newcommand{\Tc}{\mathcal{T}}
\newcommand{\Dc}{\mathcal{D}}
\newcommand{\Xc}{\mathcal{X}}
\newcommand{\Fc}{\mathcal{F}}

\newcommand{\Kh}{\widehat K}

\newcommand{\1}[1]{{\mathbf 1}{\{#1\}}}

\newcommand{\Kc}{\widehat{K}}

\let\phi=\varphi

\newcommand{\8}{{\infty}}
\newcommand{\Ber}{{B}}
\newcommand{\nn}{\nonumber}

\newcommand{\eqlaw}{\stackrel{\text{\tiny law}}{=}}
\newcommand{\cvlaw}{\stackrel{\text{\tiny law}}{\longrightarrow}}

\newcommand{\si}{\sigma}
\newcommand{\tn}{\tau_n}

\newtheorem{theo}{Theorem}[section]
\newtheorem{lm}{Lemma}[section]

\newtheorem{prop}{Proposition}[section]
\newtheorem{cor}{Corollary}[section]
\newtheorem{rem}{Remark}[section]
\newtheorem*{Th}{Theorem A}
\allowdisplaybreaks

\newcommand{\fc}{\color{black}} 
\newcommand{\nc}{\color{black}}

\title{Constrained information transmission on Erd\"os-R\'enyi graphs}
\author{Francis Comets$^{1}$ \and Christophe Gallesco$^{2}$ \and
 Serguei~Popov$^{2}$ \and Marina Vachkovskaia$^{2}$}

\begin{document}
\bibliographystyle{plain}

\maketitle  

{\footnotesize 
\noindent $^{~1}$Universit\'e Paris Diderot -- Paris 7, 
Math\'ematiques, 
 case 7012, F--75205 Paris
Cedex 13, France
\\
\noindent e-mail:
\texttt{comets@math.univ-paris-diderot.fr}

\noindent $^{~2}$Department of Statistics, Institute of Mathematics,
 Statistics and Scientific Computation, University of Campinas --
UNICAMP, rua S\'ergio Buarque de Holanda 651,
13083--859, Campinas SP, Brazil\\
\noindent e-mails: \texttt{\{gallesco,popov,marinav\}@ime.unicamp.br}

}
\bigskip

\begin{abstract}
We model the transmission of information of a message on the Erd\"os-R\'enyi random graph with parameters $(n,p)$ and limited resources.
The vertices of the graph represent servers that may broadcast a message at random. Each server has a random emission capital that decreases by one at each emission. We examine two natural dynamics: in the first dynamics, an informed server performs all  its attempts, then checks 
at each of them if the corresponding edge is open or not; in the second dynamics the informed server knows a priori 
who are its neighbors, and it performs all its attempts on its actual neighbors in the graph.
In each case, we obtain first and second order asymptotics 
(law of large numbers and central limit theorem), when $n\to \infty$ and $p$ is fixed, for the final proportion of informed servers.
 \\[.3cm]\textbf{Keywords:} 
 information transmission, rumor, labelled trees,  Erd\"os-R\'enyi random graph
 \\[.3cm]\textbf{AMS 2000 subject classifications:}
Primary 90B30; secondary 05C81, 05C80, 60F05, 60J20, 92D30
\end{abstract}

\section{Introduction}

 Information transmission with limited resources on a general graph is a natural problem which appears in various contexts and attracts an increasing interest. Consider a finite graph, each
 vertex will be seen as a server with a finite resource (e.g., operating battery) given by an
 independent random variable~$K$. Initially a message 
 comes to one of the servers, which will recast the message to its neighbors in the graph
 as long as its battery allows. In turn, each neighbor starts to emit as soon as it receives the message, and so on, in an asynchronous
 mode. The transmission stops at a finite time because of the resource constraint. A quantity of paramount  interest  is the final number of informed servers, i.e., of servers which ever receive the message.   
 
 Rumor models deal with ignorant individual (who ignore the rumor), spreaders (who know the 
 rumor and propagates it) and stiflers (who refuse to propagate it) in a population of fixed -- but
 large -- size. 
 Two  important models, usually presented in continuous time, are well-known: 
 the Maki-Thompson model \cite{MakiThompson},
and the Daley-Kendall model \cite{DaleyKendall} for which the number of eventual knowers 
obeys a law of large numbers \cite{Sudbury} and 
is asymptotically normal \cite{Pittel} {\color{black} (recently, a large deviations principle for the Maki-Thompson model was also obtained in \cite{Leb})}. Such results extend to a larger family of processes
\cite{LebMachadoRodriguez} using weak convergence theory for Markov processes. 
Though we focus here on mean-field type models, we just mention that lattice models 
lead to different questions \cite{BertacchiZucca, ColettiRodriguezSchinazi, GalloGarciaVargasRodriguez14}. On the other hand, it is well understood that the scaling limit of mean-field models are models on Galton-Watson trees, cf. \cite{AlvesETAL, Bordenave, CDS12}.
Rumor spreading models are alike epidemics propagation models, e.g. frog models \cite{AlvesMachadoPopov02, CQR07, FontesMachadoSarkar04, HoffmanJohnsonJunge14}, and the famous SIR (Susceptible, Infected, Recovered) model which has motivated a number of research papers.
See \cite{DaleyGani} for a survey.

The analysis of random graphs has recently seen a remarkable development \cite{Bollobas, Hofstad}. Such graphs yield a natural framework for rumor spreading and epidemic dissemination with more realistic applications to human or biological world \cite{NevokeeETAL}. Then, due to lack of homogeneity, setting the threshold concept on firm grounds is already a difficult problem 
\cite{IshamHardenNevokee}, and the literature is abundant in simulation experiments but poor in rigorous results.
On the Erd\"os-R\'enyi graph, the authors in \cite{FountoulakisHP} prove that the time needed for complete transmission 
in the push protocol (a synchronous dynamics without constraints) is equivalent to that of the complete graph \cite{FriezeGrimmett}
provided that the average degree is significantly larger than $(\ln n)$. In general, it is reasonable to look for
quantitative results  from perturbations of the homogeneous case. 
\fc
From the point of view of applications, the graph may be thought of as a wireless network, the
vertices of which are battery-powered sensors with a limited energy capacity. The
reader will find in Sect.1 of  \cite{CDS12} a discussion of applications to the performance evaluation of
information transmission in wireless networks.
\nc

On the complete graph, the process can be reduced by homogeneity to a Markov chain in the quadrant with absorption on the axis,
as recalled in the forthcoming Section \ref{sec:cgres}.
For a random graph, fluctuations of the vertex degrees create 
inhomogeneities which make the above description non-Markovian and  computations intractable. This can be already seen in the simplest example, the Erd\"os-R\'enyi graph. Homogeneity  is present, not in the strict sense but 
in a statistical one, and independence is deeply rooted in its construction. From many perspectives, this random graph with
fixed  positive $p$ has been proved to be very similar to the complete graph as $n$ becomes large. 
In the present paper we show that the
information transmission process on this random graph is a bounded perturbation of that on the complete graph with appropriate resource distribution. 
We will use the above mentioned similarities 
to construct couplings between information process with constraints
on the complete graph and on the Erd\"os-R\'enyi graph. Then we control the discrepancy between the two models and its propagation as the process evolves.

\medskip

\fc
In this paper, we consider two natural dynamics of the information transmission process 
on  the Erd\"os-R\'enyi graph:
\begin{itemize}
\item 
(i) an informed server performs $K$ attempts by choosing  a server at random independently at each attempts, then checks 
for each of them if the corresponding edge is open or not; 
\item (ii) the informed server knows a priori 
who are its neighbors, and it performs all its $K$ attempts on the set of its actual neighbors in the graph.
\end{itemize}
First of all we prove the existence of a threshold:  Transmission takes place at a macroscopic level if and only if
$p \IE K >1$ in the case (i), and iff $\IE K>1$ in the case (ii). Then, with positive probability, a positive proportion of servers will be informed, 
whereas in the case of the reverse inequalities, the final number of informed servers is bounded in probability.
The value of the threshold is natural, observing that, in the first case,  attempts taking place on closed edges are lost, so that 
the effective number of attempts is close (as $N$ increases)  to a random sum
\begin{equation}
\label{eq:Kh}
\widehat K=\sum_{1\leq k \leq K} B_k
\end{equation}
with $(B_k, k \geq 1)$ i.i.d.~Bernoulli with parameter $p$. 

Our main results, Theorems \ref{th:LLN}, \ref{th:TLC}, \ref{th:LLN2} and 
\ref{th:TLC2} below,  are the laws of large numbers and the central limit theorems for the number of informed servers  with explicit values of the limits in each case.
Our approach is to show that, in the limit $n \to \8$ with a fixed $p \in (0,1]$, 
the information transmission process on the Erd\"os-R\'enyi graph is shown to be a bounded perturbation of the 
process on the complete graph with a suitable resource law. 
Then, the first and second order asymptotics, obtained by explicit computations on the complete graph in \cite{Ma3, CDS12}, still hold on the random graph. 
\nc

An important property of the model is abelianity, 
\fc e.g. see Proposition \ref{prop:2cg}.
\nc
 We can change the order in which emitters are taken without changing 
\fc   the law of the final state of the process, 
and construct an efficient  
coupling of the processes on the two graphs. This property also implies that assuming the servers emit in a burst does not change the final result ({\color{black} a nice feature of the burst emission assumption is that it reveals a branching structure}).
The two dynamics we consider here are simple and reasonable protocols, but we don't make any attempt  for generality in this paper.
\nc
We will use an exploration process which allows to reveal at each step, only the necessary part of the graph in order to preserve randomness and stationarity in the subsequent steps.

{\bf Outline of the paper}: In Section 2 we define the model, recall useful results for the complete graph, and state our main results. Then, labeled trees are introduced with a view towards our constructions.  Section 4 contains the proofs in the case of the first dynamics (i), and the last section deals with dynamics (ii).

\section{Model and results}

We start to recall some results for the information transmission process on the complete graph.



\subsection{Known asymptotics in the case of the complete graph}	\label{sec:cgres}

When any server is connected to any other one, the communication network is the complete graph on ${\cal N}=\{1,\ldots,n\}$.
{\color{black} We consider here discrete time and} we scale the time so that there is exactly one emission per time unit. 
Then, the  information process  can be fully described by  the number   $N_n(s)$ of informed servers at time $s$ and the   
number $S_n(s)$ of available emission attempts
\fc
 (see \cite{CDS12} for the formal definition).
\nc
 Precisely, for the information process with resource $K$ on the complete graph,  the pair  
 $(S_n(s),N_n(s))_{s=0,1,\ldots}$  on $\Z_+ \times [1,n]$ is  a Markov chain 
with transitions 
\begin{equation} \label{eq:CDSdyn}
\left\{
\begin{array}{rcl}
 \IP\Big( S_n(s \! + \! 1)=S_n(s) \! - \! 1,N_n(s \! + \! 1)=N_n(s) \mid \Fc_s\Big) &  =&\frac{N_n(s)}{n},   \\
 \IP\Big( S_n(s \! + \! 1)=S_n(s) \! + \! k \! - \! 1,N_n(s \! + \! 1)=N_n(s) \! + \! 1 \mid \Fc_s\Big) &  =&\left(1-\frac{N_n(s)}{n}\right) \IP(K=k),\
\end{array}
\right.
\end{equation}
for $k \geq 0,$ with $\Fc_s$ the $\sigma$-field generated by $S_n(\cdot)$ and $N_n(\cdot)$
on $[0,s]$. 
\fc
The  transition probabilities are easily understood by interpreting  what can occur at a given step:
\nc
On the first line {\color{black} of (\ref{eq:CDSdyn})} the emission takes place towards a previously informed target,  though in the second one
the target yields its own resource (a fresh r.v.~$K$). 
The chain is absorbed in  the vertical semi-axis, at the finite time ${\mathfrak T}_n = \inf\{s: S_n(s)=0\}$.
\medskip

In this section we recall some results from \cite{CDS12} (and of \cite{Ma3} for constant $K=2$) on the first and second order 
asymptotics of $N_n({\mathfrak T}_n )=N_n(\8)$.
Let $q \in [0,1)$ be the largest root of 
\begin{equation}
  \label{eq:p(theta)}
  q\; \IE K +\ln (1-q) =0.
\end{equation}
Then, $0<q<1 $ for $\IE K>1$ and $q=0$ if $\IE K \leq 1$.
\begin{Th}[\cite{CDS12}, Theorems 2.2 and 2.3]{}\label{th:cdsLLNTCL}
(i)  Assume $\IE K \in (0,\8)$. 
 Then, as $n \to \8$, 
$$
\frac{1}{n}
    N_n({\mathfrak T}_n )
\cvlaw 
q \times Ber(\sigma^{GW}) 
$$
with $Ber(\sigma)$ a Bernoulli variable with parameter $\sigma$, and 
$\sigma^{GW}$ is the largest solution $\sigma \in [0,1]$ of 
\begin{equation}
\label{eq:probasurvivalgw}
1-\sigma= \IE\left[(1-\sigma)^K\right],
\end{equation}
i.e. the survival probability of a Galton-Watson process with reproduction law $K$. 
\\

(ii)
Assume $\IE K >1$ and $\IE K^2 < \8$. 
Denote by  $\si_K^2$ the variance of $K$ and fix some $\eps$  with $0< \eps <- \ln (1-q)$. As $n \to \8$, we have the convergence in law, conditionally on $\{{\mathfrak T}_n \geq \eps n\}$,
$$
n^{-1/2} \big( N_n({\mathfrak T}_n) - n q \big) \cvlaw {\mathcal N}(0, \si_q^2),
$$
with ${\mathcal N}(0, \si^2)$ a centered Gaussian with variance $\si^2$, and
\begin{equation}
\label{eq:varTCLq}
\si_q^2=\frac{ q\si_K^2 (1-q)^2+  q(1-q) + (1-q)^2 \ln(1-q)}{[(1-q)\IE K-1]^{2}}
\fc
>0.
\nc
\end{equation}
\end{Th}

We now state the main results of this paper, i.e. when the connection network is the Erd\"os-R\'enyi graph
$G(n,p)$.  Now, a server starting to emit, instantaneously exhausts its $K$ emissions in a burst. 
The time unit corresponds to complete exhaustion for an emitter.      
Let $N_n^{er}(t)$ be the number of informed servers at time $t$.

\subsection{
First mode of transmission on the  Erd\"os-R\'enyi graph} \label{sec:res1}

First of all, the  Erd\"os-R\'enyi graph $G(n,p)$ is sampled on the vertex set $\cal N$ (each unoriented edge 
is kept independently 
with probability $p$ or removed with probability $1-p$), and one vertex is selected
as the first informed server. Then, at each integer time, an   informed server which is not yet exhausted is selected
to emit its $K$ attempts in a burst. For each attempt a target in $\cal N$ is selected (in the full population including the emitter). If the target is already informed or if the corresponding edge is not in the graph, the attempt is lost.
Otherwise, the target becomes informed. After all attempts are checked, the emitter is turned to exhausted and the 
time is increased by one unit. The transmission ends at a finite time $ \tau_n^{er}$, 
\fc which is the first time when all informed servers are exhausted.
\nc
\medskip

Note that, because of the burst emission here,  the {\it time scale is different} from Section \ref{sec:cgres} with one emission at a time. 
\fc With $N_n^{er}(t)$ the number of informed servers at time $t$, 
\nc
we are interested in the asymptotics of 
$$\tau_n^{er} =
N_n^{er}( \tau_n^{er} ) = N_n^{er}(\8).$$
\fc The first equality holds since it takes one time unit to exhaust an informed server, and the last one holds since the process stops at $\tau_n^{er} $.
\nc

We will encounter the above quantities when $K$ is replaced by $\widehat K$ from (\ref{eq:Kh}), that we will denote using the same symbol with a hat: 
\fc In particular,
\nc
 $\widehat \theta =0=\widehat \sigma^{GW}$ if $p \IE K \leq 1$, and for $p \IE K >1 , \widehat q \in (0,1)$ is the positive root of 
\begin{equation} \label{eq:hatq}
\widehat q p \IE K  + \ln (1-{\widehat q})=0,
\end{equation}
and $\widehat \sigma^{GW}$ is the positive root of
\begin{equation*}
1-\widehat \sigma^{GW}= \IE\left[(1-\widehat \sigma^{GW})^{\widehat K}\right]=
\IE\left[(1-p \widehat \sigma^{GW})^{K}\right],
\end{equation*}
\fc   that is, equation (\ref{eq:probasurvivalgw}) with hats.
\nc
\begin{theo} \label{th:LLN}
Assume $\IE K^2<\8$. Then,  
$$
\frac{  \tau_n^{er} }{n}
\cvlaw   \widehat q \times
Ber(\widehat \sigma^{GW})  
\;.
$$
\end{theo}
The interesting case is of course when $\IE \widehat K=p \IE K >1$ to have $\widehat \sigma^{GW}>0$. In this case, let also
\begin{eqnarray*}
	\widehat \si_{{\color{black}\hat{q}}}^2=\frac{ \widehat q\si_{\widehat K}^2 (1-\widehat q)^2+  \widehat q(1-\widehat q) + (1-\widehat q)^2 \ln(1-\widehat q)}{[(1-\widehat q)p \IE K-1]^{2}}{\color{black}>0},
\end{eqnarray*}
with $\si_{\widehat K}^2= p(1-p) \IE K + p^2 \si_K^2 $ the variance of $\widehat K$.

\begin{theo} \label{th:TLC}
Assume $\IE K^2<\8$ and $p \IE K >1$. Fix $\eps \in (0,\widehat q)$. Then,
conditionally on $\{\tau^{er}_n > n \eps\}$, we have convergence in law:
$$
\frac{
 \tau_n^{er} 
- n \widehat q
}{\sqrt n}
\cvlaw 
{\mathcal N}(0, \widehat \sigma_{{\color{black}\hat{q}}}^2)\;.
$$
\end{theo}

\subsection{Main results for the second mode of transmission} \label{sec:res2}

Again, we start by sampling the  Erd\"os-R\'enyi graph $G(n,p)$  and one vertex 
as the first informed server. Then, at each integer time, an  informed server which is not yet exhausted is selected
to emit its $K$ attempts in a burst, each attempt being towards a random  target uniformly distributed among the 
neighbors in the graph. 
(If a site has no neighbours, it wastes its resource without result, and after that the process 
continues.)
If the target is already informed  the attempt is lost, but 
otherwise the target becomes informed. After all attempts are checked, the emitter is turned to exhausted and the 
time is increased by one unit. The transmission ends at some finite time $ \bar \tau_n^{er}$,
with $\bar \tau_n^{er}$ informed servers.
In the following theorems $q$, $\sigma^{GW}$ and $\sigma_q$ are from Section \ref{sec:cgres}.

\begin{theo} \label{th:LLN2}
Assume $\IE K^2<\8$. 
$$
\frac{\bar \tau^{er}_n}{n}
\cvlaw 
q \times Ber( \sigma^{GW})  .
$$
\end{theo}
\begin{theo} \label{th:TLC2}
Assume $\IE K^2<\8$ and $\IE K >1$. Fix $\eps \in (0,q)$. Then,
conditionally on $\{\bar \tau^{er}_n > n \eps\}$, we have convergence in law: 
$$
\frac{
\bar \tau^{er}_n
- n q
}{\sqrt n}
\cvlaw 
{\mathcal N}(0, \sigma_{q}^2).
$$
\end{theo}

\subsection{Strategy of the proofs}

\fc
We use the known results about the information process on the complete graph to derive results on the Erd\"os-R\'enyi graph. 

We show that case (i) is similar to the complete graph with $\hat K$ attempts. 
The difference is that in the 
latter model, the Bernoulli random variables in (\ref{eq:Kh})
(indicating the presence of the relevant edges) 
are regenerated
independently at each attempt to transmit, though in the former the state of an edge is determined at its first appearance. A coupling argument is made to show that in fact this makes little difference to the final number of vertices receiving the information. 
 In case (ii) we keep track of which edges are in a known state and the key argument is that with high probability only $o(n)$ edges out of a vertex will ever be in a known state.
Hence the argument is to show that, most likely, there will be $O_P(1)$ transmissions in which there is a discrepancy between the models. To take care of the consequences of discrepancies, we delay them until the end of the process -- taking advantage of irrelevance of the order of transmission. Finally we show that these few extra transmissions make little difference to the final proportion of vertices receiving the information. 
\nc

\section{Construction from labelled trees}

\subsection{Labelled trees}

Let $\Wc=\cup_{m\geq 0}\N^m$ be the set of all finite words on the alphabet $\N=\{1, 2,\dots\}$. By convention $\N^0=\{ \Par \}$ 
contains only one element which can be interpreted as the empty word and which, in our formalism, will be the root of the tree. An element of $\Wc$ different from $\Par$ is thus a $m$-uple $u=(i_1,\dots,i_m)$ which, to simplify, will be denoted by $u=i_1\dots i_m$. The length of $u$ denoted by $|u|$ equals $m$ (with $|\Par|=0$). If $j\in \N$, we denote by $uj$ the element $i_1\dots i_mj$. The elements of the form $uj$ are interpreted as the descendants of $u$.  We will use the following total order relation on $\Wc$, we write $w\leq w'$ if: $|w|<|w'|$, or $|w|=|w'|$ and $w\leq_{lex} w'$ in the lexicographical order. 

A {\it rooted tree} $\Tc$  is an undirected simple connected graph without cycles and with a distinguished vertex. A {\em labelled tree} $(\Tc,L)$ is a rooted tree $\Tc$ 
equipped with a label mapping $L$ from
$\Tc$ to some set $\Nc$. Labelled trees we will consider are connected subsets of $\Wc$
containing $\Par$. The label set  $\Nc=\{1,2,\ldots,n\}$ encodes the set of servers.
For the sake of brevity, we  use the short notation $Lv=L(v)$, that the reader will distinguish from concatenation.

\subsection{Construction and coupling} \label{sec:constr}

Let $n\geq 2$ and $\Nc=\{1,\dots,n\}$. On a suitable probability space $(\Omega, \mathcal{F}, \IP)$, we define the following independent random elements:
\begin{itemize}
\item[(i)] $(K_i)_{i\in \Nc}$ are i.i.d.\ non-negative integer random variables;
\item[(ii)] $I_0$ is a uniform random variable on $\Nc$;
\item[(iii)] $(I^i_k)_{(i,k)\in \Nc\times \N}$ are independent uniform random variables on $\Nc$;
\item[(iv)] $(\Ber^i_k)_{(i,k)\in \Nc\times \N}$ are independent Bernoulli random variables of parameter $p$.
\end{itemize}
\medskip

In the next sections, we construct couplings between different labelled trees using the above random elements. The different trees $\Tc(\8)= \lim_{t \nearrow \8}\Tc(t)$ are limits of some 
sequence, they
will be constructed dynamically discovering step by step their nodes and their labels.
At each step $t$, the labels of the tree $\Tc(t)$ are all different.
\fc 
They represent the 
{\em informed}  servers at time~$t$ and 
will be partitioned
\nc
 into two (as in Section \ref{sec:er}) or three (as in Section \ref{sec:2frigos} and the end of Section \ref{sec:cg}) subsets:
 \begin{equation}
\label{eq:partition}
\Ic(t) = \Ac(t) \cup \Ec(t) \cup \Dc(t), \qquad \Ac(t), \Ec(t), \Dc(t) \; {\rm disjoint},
\end{equation}
where $\Ec(t)$ encodes the exhausted servers (those which have already used their resource),
$\Ac(t)$ encodes the active servers (those which are waiting to use 
\fc their resource and  ready to transmit),
\nc
 and $\Dc(t)$ 
\fc
encodes the set of delayed servers
(those which have not started to transmit but are temporarily delayed). 
\nc
In Section \ref{sec:er}, $\Dc(t)$ is empty.

\begin{rem}
The sets in (\ref{eq:partition}) and the mapping $L$ depend on the number $n$ of servers. In general, for the sake of simplicity,  we do not indicate explicitly the dependence in the notations.
\end{rem}

In the next section, we construct  transmission processes on the complete graph with law $\Kc$ and on Erd\"os--R\'enyi random graph, using the above elements,  
thus we have a coupling between these processes, allowing to transfer results from one to the other. 
\fc The reader may wonder, here or below, why we introduce so many independent r.v.'s in the construction,
since a given edge is decided to be open or closed in the  Erd\"os-R\'enyi graph only once,
namely at its first appearance. The reason is that coupling the process with one on the complete graph requires to decide each edge more than once
(cf.~ Sections \ref{sec:er} and \ref{sec:cg}).
\nc
\section{First mode of emission: Burst emission}

We model the transmission of a message on the Erd\"os-R\'enyi random graph of parameter $p$. Each vertex $i$ of the graph is a server with resource $K_i$. Initially, one vertex  receives a message and tries to send it to its neighbors: first it choses, uniformly among all the servers, one server (the target) to which it will try to send the message.  If the edge between these two servers is present, then the information is transmitted and otherwise it is not. The emitting server repeats this operation until it has 
exhausted its own resource $K_i$. If the edge is present and the target server already knows the information then the emitter just loses one  resource unit. When the emitter has exhausted its resource, we pick a new server among the informed ones and it starts to emit according to  the same procedure. The process stops when all the informed servers have exhausted their resources.

\subsection{Erd\"os-R\'enyi graph} \label{sec:er}


\begin{figure}[!htb]
\begin{center}
\includegraphics[scale= 1.0]{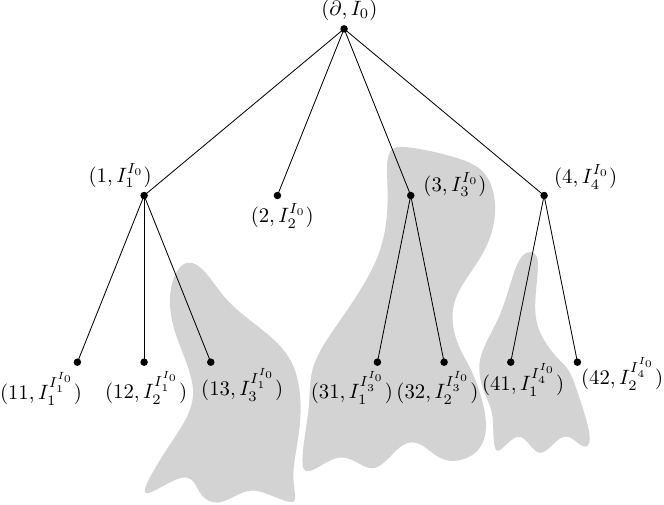}
\caption{Two first generations of 
a labeled random  tree, with $K_{I_0}=4,$ $ K_{I_1^{I_0}}=~3, 
K_{I_2^{I_0}}=0,$  $K_{I_3^{I_0}}=2,$ and $ K_{I_4^{I_0}}=2$ is represented by the line segments.
The other variables used in the construction are: $n=38; I_0=23; I^{23}_1=12, I^{23}_2=7, I^{23}_3=12, I^{23}_4=30, I^{12}_1=18, I^{12}_2=I^{12}_3=2, I^{30}_2=27; 
B^{23}_k=1 (1\leq k \leq 4), B^{12}_1=B^{12}_2=1, B^{12}_3=0, B^{30}_1=0, B^{30}_2=1$. Some parts of the tree, included in the shaded parts,  
were pruned out, because a label has been already discovered or because an edge is closed in the Erd\"os-R\'enyi graph. 
}
\label{fig1}
\end{center}
\end{figure}

We construct dynamically the random labelled tree $\Tc^{er}(t)$ in the following way. 
At $t=0$, using $I_0$, we discover the label of the root $\Par$ and we set $L(\Par)=I_0$. In the rest of this section, we abbreviate for clarity $L^{er}$ by $L$.
Denoting by $i$ the value of $I_0$ for short notations, 
we consider a realization of $K_i$, $\Ber^i_1,\dots, \Ber^i_{K_i}$ and $I^i_1,\dots,I^i_{K_i}$. To each descendant $l$, $1\leq l \leq K_i$, of the root we associate 
the label $I^i_{l}$. Initially, we define
\begin{align}
\nonumber
X(0)&=\Par,\\
\Xc(0)&=\{w\in \Wc : |w|=1, 1\leq w\leq K_i\},\nonumber\\
\Xc^{er}(0)&=\{w\in \Wc : |w|=1, 1\leq w\leq K_i, \Ber^i_w=1\},\nonumber\\
L(w) &= I^i_w, \qquad w \in \Xc(0),\nonumber\\
\Ac^{er}(0)&=\{w\in \Xc^{er}(0) :  L(w)\neq i, L(w)\neq L(w'), w'<w, w'\in
 \Xc(0)\},
\nonumber\\
\Ic^{er}(0)&=\{\Par\}\cup \Ac^{er}(0).  \label{eq:er0}
\end{align}
With the process $(X(t), \Ic^{er}(t), 
\Ac^{er}(t))$ at time $t$ and $L$ defined on $\Tc^{er}(t)$, the value at the next step $t+1$ is defined by:
\begin{itemize}
\item If $\Ac^{er}(t)$ is non empty, we let $X(t+1)$ be its first element in the total order $\leq$,
\[
X(t+1)=\inf\{w\in \Ac^{er}(t)\},\phantom{*}\mbox{denoted by $v$}
\]
and we consider a realization of $K_{Lv}$, $\Ber^{Lv}_1,\dots, \Ber^{Lv}_{K_{Lv}}$ and $I^{Lv}_1,\dots,I^{Lv}_{K_{Lv}}$. To each descendant $vl$, $1\leq l \leq K_{Lv}$ we associate the label $I^{Lv}_{l}$. Then, we update the sets of vertices
\begin{align}
\Xc(t+1)&=\{vl \in \Wc : 1\leq l \leq K_{Lv}\},\nonumber\\
L(w) &= I^{Lv}_l, \qquad w=vl \in \Xc(t+1),\nonumber\\
\Xc^{er}(t+1)&=\{vl \in \Wc : 1\leq l \leq K_{Lv}, \Ber^{Lv}_l=1\},\nonumber\\
\Ac^{er}(t+1)&=(\Ac^{er}(t)\setminus \{v\})\cup\{w\in \Xc^{er}(t+1) :  L(w)\notin L(\Ic^{er}(t)), \nonumber\\
&\qquad \qquad \qquad L(w)\neq L(w'), w'<w, w'\in \Xc(t+1)\} ,\nonumber\\
\Ic^{er}(t+1)&=\Ic^{er}(t)\cup \Ac^{er}(t+1).
\label{eq:ert}
\end{align}

\item If $\Ac^{er}(t)$ is empty, we set $\tau_n^{er}=t$ and the construction is stopped.
\end{itemize}

At each step of the construction, 
we set $\Ec^{er}(t+1)=\Ec^{er}(t) \cup \{X(t+1)\}$ starting from $\Ec^{er}(0)=\{\Par\}$,
and $\Dc^{er}(t)\equiv \emptyset$.
Hence $v=X(t+1)$ is moved from active to exhausted at time $t+1$, and 
the   partition (\ref{eq:partition}) reduces to two subsets
in the case of the Erd\"os-R\'enyi graph.   \qed
\medskip

The construction is illustrated by Figure \ref{fig1}.

\begin{rem} \label{rem:L}
\fc  (i) Note that the definitions of active servers in (\ref{eq:er0}) and  (\ref{eq:ert}) 
  require that the label has not appeared before. Indeed the first Bernoulli variable determines the status of the edge.
\nc

(ii)
Note that $L^{er}$ has been defined in this process on a larger tree than needed. In fact, we consider its restriction  to $\Tc^{er}(\8)=\Tc^{er}(\tau_n^{er})$, which is indeed injective.
This procedure of restriction is needed in all the subsequent constructions.
\end{rem}

It is not completely obvious that this construction corresponds to the description given at the beginning of Section  \ref{sec:res1}. However, this is the case, as we show now. Here, an edge is open or closed according to the Bernoulli variable used on the first
appearance of the edge in the construction. Here is a formal definition.
Denote by $\cal E$ the set of unoriented edges on $\cal N$, i.e. the set of 
$e = \langle i,j \rangle, i, j \in {\cal N}$ (allowing self-edge). We say that the edge $e$ has appeared in the construction if 
there is some $t \;(0 \leq t \leq \tau_n^{er}-1)$ and some $\ell \leq K_{L[X(t)]}$ such that 
$$e= \langle L[X(t)], L[X(t)\ell] \rangle.$$
We denote by $t(e)$ and $\ell(e)$ the smallest (in the lexicographic order) 
$t$ and $\ell$ with the above property, by $Ap(t)=\{e: t(e)=t\}$ the set of edges which have appeared at time $t$ and 
$Ap=\cup_{t \geq 0} Ap(t)$.
With an additional i.i.d. Bernoulli($p$) family $(\Ber^i_{-k})_{i \in {\cal N}, k \geq 1}$
independent of the variables in (i--iv), define
$$
B(e) = 
\left\{
\begin{array}{lll}
B^i_k  & {\rm if}  &  e \in Ap, i=L[X(t(e))], k=\ell(e), \\
B^i_{-j} & {\rm if}  &  e \notin Ap, e=\langle i,j \rangle, i \leq j. 
\end{array}
\right.
$$
\begin{prop} \label{prop:er=er}
The family $(B(e), e \in {\cal E})$ is i.i.d.~Bernoulli($p$), and is independent of the family
$(I^i_k)_{i \in {\cal N}, k \geq 1}$,$ (K_i)_{i \in {\cal N}}$, $I_0$. 
\end{prop}
\noindent
The proposition shows that the above construction coincides with the description of the information transmission process on the 
Erd\"os-R\'enyi graph as given in  the beginning of Section~\ref{sec:res1}, in which the random graph is defined by the $B(e)$'s, and the dynamics uses the variables $I_0, I^i_\cdot, K_i$. Though the proof is standard, we give it for completeness.
\medskip

$\Box$ For bounded measurable functions $f_e, g_i, h$ defined on the appropriate spaces, we compute
\begin{eqnarray*}
\IE  \prod_{e \in {\cal E}} f_e(B(e)) \times  \prod_{i \in {\cal N}} g_i(K_i, I^i_\cdot) \times  h(I_0)  
=\qquad \qquad \qquad \qquad \qquad \qquad \qquad \\
\sum_{T, A(\cdot)} 
\IE  \prod_{e \in {\cal E}} f_e(B(e)) \times  \prod_{i \in {\cal N}} g_i(K_i, I^i_\cdot) \times  h(I_0) 
{\mathbf 1}_{\{\tau_n^{er}=T, Ap(\cdot)=A(\cdot)\}}
\end{eqnarray*}
where $A(\cdot)=(A(t), t=0,\ldots T-1)$ ranges over the collection of $T$ disjoints subsets of $\cal E$.
With $A=\cup_t A(t)$, by independence of  $(\Ber^i_{-k})_{i,k}$
with the other variables, the expectation in the last term is equal to
\begin{eqnarray} \label{eq:Kbg}
 \Big( \prod_{e \notin A} \IE  f_e(B(e)) \Big) \times \Big( \IE  \prod_{e \in A}   f_e(B(e)) \times  
 \prod_{i \in {\cal N}} g_i(K_i, I^i_\cdot) \times  h(I_0) 
{\mathbf 1}_{\{\tau_n^{er}=T, Ap(\cdot)=A\}} \Big).
\end{eqnarray}
Let us write the last expectation as
\begin{eqnarray*}
   \IE \Big(\prod_{e \in A(T-1)}  f_e\big(   B^{L[X(t(e))]}_{ \ell(e)}\big) \times  
   \prod_{e \in A \setminus A(T-1)}  f_e(B(e))  \times
 \prod_{i \in {\cal N}} g_i(K_i, I^i_\cdot) \times  h(I_0) 
{\mathbf 1}_{\{\tau_n^{er}=T, Ap(\cdot)=A\}} \Big)
\qquad
\\
\stackrel{\rm indep.}{=}
 \Big( \prod_{e \in A(T-1)}  \IE f_e(   B(e))\Big) \times \Big( \IE
   \prod_{e \in A \setminus A(T-1)}  f_e(B(e))  \times
 \prod_{i \in {\cal N}} g_i(K_i, I^i_\cdot) \times  h(I_0) 
{\mathbf 1}_{\{\tau_n^{er}=T, Ap(\cdot)=A\} }\Big)\\
\stackrel{\rm iterating}{=} 
\Big( \prod_{e \in A}  \IE f_e(   B(e))\Big) \times \Big( \IE
 \prod_{i \in {\cal N}} g_i(K_i, I^i_\cdot) \times  h(I_0) 
{\mathbf 1}_{\{\tau_n^{er}=T, Ap(\cdot)=A\}} \Big).
\qquad\qquad\qquad\qquad\qquad
\end{eqnarray*}
The first factor in the last  term complements the one in (\ref{eq:Kbg}), and we finally get
\begin{align*}
\lefteqn{\IE  \prod_{e \in {\cal E}} f_e(B(e)) \times  \prod_{i \in {\cal N}} g_i(K_i, I^i_\cdot) \times  h(I_0)}\phantom{**************}\nonumber\\
&=
\prod_{e \in {\cal E}}  \IE f_e(   B(e))  \Big( \sum_{T, A(\cdot)}   \IE
 \prod_{i \in {\cal N}} g_i(K_i, I^i_\cdot) \times  h(I_0)  
{\mathbf 1}_{\{\tau_n^{er}=T, Ap(\cdot)=A\}} \Big) \\
&=
 \prod_{e \in {\cal E}}  \IE f_e(   B(e)) \times \Big( \IE
 \prod_{i \in {\cal N}} g_i(K_i, I^i_\cdot) \times  h(I_0)  \Big),
\end{align*}
by summing over $T, A(\cdot)$. This proves the proposition. \qed


\subsection{Complete graph} \label{sec:cg}

{\color{black} With the above ingredients, we start to construct the information transmission model on the complete graph according to two different dynamics with distribution $\Kh$. The first construction is simple and natural (and is close to that performed in Section 2.1 of \cite{CDS12}),} but the second one makes a useful
coupling with the transmission model on the Erd\"os-R\'enyi graph.

\medskip

{\bf Sequential construction.}
We construct dynamically the random tree $\Tc^{cg,s}(t)$  together with  the label mapping $L^{cg,s}$ that we abbreviate for clarity by $L=L^{cg,s}$ in the construction. The exploration vertex 
$X(t)$ we introduce below, also depends on the dynamics, $X=X^{cg,s}$, but we omit the superscript for the same reason. 

At $t=0$, using $I_0$, we discover the label of the root, $L(\Par)=I_0$. Suppose that $I_0=i$, and consider the realization of $K_i$, $\Ber^i_1,\dots, \Ber^i_{K_i}$ and $I^i_1,\dots,I^i_{K_i}$. The label of each descendant $l$, $1\leq l \leq K_i$, of the root 
is $I^i_{l}$. Initially, we suppose 
\begin{align}
\nonumber
X(0)&=\Par, \\
\Xc^{cg,s}(0)&=\{w\in \Wc : |w|=1, 1\leq w\leq K_i, \Ber^i_w=1\},\nonumber\\
L(w) &= I^i_w, \qquad w \in \Xc^{cg,s}(0),\nonumber\\
\Ac^{cg,s}(0)&=\{w\in \Xc^{cg,s}(0) :  L(w) \neq i, L(w)\neq L(w'), w'<w, w'\in \Xc^{cg,s}(0)\},\nonumber\\
\Ic^{cg,s}(0)&=\{\Par\}\cup \Ac^{cg,s}(0). \label{eq:cgs0}
\end{align}
With the process $(\Ic^{cg,s}(t),\Xc^{cg,s}(t), \Ac^{cg,s}(t))$ at time $t$, its value at the next step $t+1$ is defined by:
\begin{itemize}
\item If $\Ac^{cg,s}(t)$ is non empty, we let $X(t+1)$ be its first element,
\[
X(t+1)=\inf\{w\in \Ac^{cg,s}(t)\},\phantom{*}\mbox{denoted by $v$}
\] 
and we consider a realization of $K_{Lv}$, $\Ber^{Lv}_1,\dots, \Ber^{Lv}_{K_{Lv}}$ and $I^{Lv}_1,\dots,I^{Lv}_{K_{Lv}}$. To each descendant $vl$, $1\leq l \leq K_{Lv}$ we associate the label $I^{Lv}_{l}$. Then, we update the sets:
\begin{align}
\Xc^{cg,s}(t+1)&=\{vl \in \Wc : 1\leq l \leq K_{Lv}, \Ber^{Lv}_l=1\},\nonumber\\
L(w) &= I^{Lv}_l, \qquad w=vl \in \Xc^{cg,s}(t+1),\nonumber\\
\Ac^{cg,s}(t+1)&=(\Ac^{cg,s}(t)\setminus \{v\})\cup\Big\{w\in \Xc^{cg,s}(t\!+\!1) :  L(w)\notin L(\Ic^{cg,s}(t)), \nonumber\\
&\phantom{**************}L(w)\neq L(w'), w'<w, w'\in \Xc^{cg,s}(t\!+\!1)\Big\},\nonumber\\
\Ic^{cg,s}(t+1)&=\Ic^{cg,s}(t)\cup \Ac^{cg,s}(t+1). \label{eq:cgst}
\end{align}

\item If $\Ac^{cg,s}(t)$ is empty, we set $\tau_n^{cg,s}=t$ and the transmission stops.
\end{itemize}

At each step of the construction, 
we set $\Ec^{cg,s}(t+1)=\Ec^{cg,s}(t) \cup \{X(t+1)\}$ starting from $\Ec^{cg,s}(0)=\{\Par\}$,
and $\Dc^{cg,s}(t)\equiv \emptyset$, so that the   partition (\ref{eq:partition}) reduces again to two subsets. 
\qed

\begin{rem} 
{\color{black} Observe that the condition ``$L(w)\neq L(w')$, $w'<w$" in the definition of $\Ac^{cg,s}(t)$ avoids counting twice the same label.} We also emphasize that \eqref{eq:cgs0}--\eqref{eq:cgst} and  \eqref{eq:er0}--\eqref{eq:ert} differ by the set 
$\Xc^\cdot(t)$ for $w'$ in the next-to-last line of the formulae. In the Erd\"os-R\'enyi construction the same Bernoulli
variable is in force each time an edge is used, whereas a fresh Bernoulli
variable is needed on the complete graph.
\end{rem}

\medskip

In the next proposition we show that this construction  yields the information transmission model on the complete graph from \cite{CDS12} with distribution $\Kh$.  This fact is necessary in order to use known results on the complete graph. 
\fc The time scales differently in the two constructions. 
\nc
To relate them, we introduce, for all $i\in \Nc$, 
\begin{equation}\nonumber
\widehat K_i = \sum_{k=1}^{ K_i} B^i_k,
\end{equation}
and, for $t=0,1,\ldots \tau_n^{cg,s}$,  
$\widehat R=\widehat R^{cg,s}_n$ by
\begin{equation}\nonumber
\widehat R(t)= \sum_{r=0}^{t-1} \widehat K_{LX(r)},\qquad \widehat R(0)=0,
\end{equation}
where we recall that $X=X^{cg,s}$ and $L=L^{cg,s}$. We also define $\widehat N_n^{cg,s}, \widehat S_n^{cg,s}$ starting from  the initial configurations
$\widehat N_n^{cg,s}(0)=1, \widehat S_n^{cg,s}(0)= \hat{K}_{L\Par}$,  with the following evolution.
\begin{itemize}
\item For $s \in ]\widehat R(t),\widehat R(t+1)]$ with $t \in [0,\tau_n^{cg,s}[$, we consider the smallest integer 
$\ell = \ell_s\in [1, K_{LX(t)}]$ such that $ \sum_{k=1}^{ \ell} B^{LX(t)}_k=s-\widehat R(t)$, and we define
\begin{eqnarray}
\label{eq:dodobiento}
\widehat N_n^{cg,s}(s) \!\!\!\!
&=&  \!\!\!\!
\widehat N_n^{cg,s}(s-1)+ 
\1{X(t)\ell  \in \Ac^{cg,s}(t)},\\
\nn 
\widehat S_n^{cg,s}(s)  \!\!\!\!
&=&  \!\!\!\!
\widehat  S_n^{cg,s}(s-1)+  \1{X(t) \ell \in \Ac^{cg,s}(t)}
\times  \widehat K_{L\big(X(t)\ell\big)}-1 ,
\end{eqnarray}
where $X(t) \ell$ denotes by concatenation a direct child of $X(t)$ in the tree.
We  check from \eqref{eq:NetSaRt} below that 
\begin{equation}
\label{eq:flb}
\widehat R( \tau_n^{cg,s})= \inf\{ s \geq 0: \widehat  S_n^{cg,s}(s)=0\}.
\end{equation}
\item After that time the process stops: $\widehat N_n^{cg,s}(s)=\widehat N_n^{cg,s}(\widehat R(\tau_n^{cg,s}))$ and $\widehat S_n^{cg,s}(s)=0$
for $s \geq \widehat R(\tau_n^{cg,s})$.
\end{itemize}

Note, for further use, that by summing \eqref{eq:dodobiento}, we find for all $t$,
\begin{eqnarray}
\nn 
\widehat  N_n^{cg,s}(\widehat R(t))&=&{\rm card} \;\Tc^{cg,s}(t), \qquad t \geq 0,\\
\label{eq:NetSaRt}  
\widehat S_n^{cg,s}(\widehat R(t))&=& \sum_{u \in \Ac^{cg,s}(t)}  \widehat K_{Lu} , \qquad t \leq \tau_n^{cg,s}.
\end{eqnarray}

The process $(\widehat S_n^{cg,s}, \widehat N_n^{cg,s}) (\cdot)$ is the one considered in \cite{CDS12}, i.e. we recover the dynamical definition of the information
transmission process on the complete graph:

\begin{prop} \label{prop:cgs}
$ \left(
\widehat S_n^{cg,s} (s) ,
\widehat N_n^{cg,s}  (s)
\right)_{s\geq 0} $
 is a  Markov chain  with transitions given by  (\ref{eq:CDSdyn})  and  resource variable
$\widehat K$,  stopped when the first coordinate vanishes. 
Moreover, we have equality in law 
$$(\widehat R( \tau_n^{cg,s}),  {\rm card} \;\Tc^{cg,s}(\8))
 \eqlaw ({\mathfrak T}_n,  N_n (\8)).$$
\end{prop}
$\Box$ From the independence of the random elements (i)--(iv) in Section \ref{sec:constr}, it  is a standard exercise to check it is a Markov chain, and from (\ref{eq:dodobiento}) that the transition probability 
 is given by
 \begin{eqnarray*}
 (\sigma, \nu) \longrightarrow 
 (\sigma + B \widehat K -1, \nu +B),
 \end{eqnarray*}
where $\widehat K$ and $B$ are independent variables with law (\ref{eq:Kh}) and Bernoulli with parameter $1-\nu/n$.
Together with the absorption rule at the random time defined by (\ref{eq:flb}), this proves the first claim.
The other one then follows from (\ref{eq:NetSaRt}) and \eqref{eq:flb}.
\qed

\medskip

\fc  Propositions \ref{prop:er=er} and \ref{prop:cgs}  mean that we have a coupling between the information transmission process on the two graphs.
\nc
The main problem with it is that after the first discrepancy between $\Ic^{er}(t)$ and $\Ic^{cg,s}$ 
occurs, the two constructions diverge and we loose  track of the differences. For instance, 
we have $\Ac^{er}(t) \subset \Ac^{cg,s}(t)$ for $t=0$, 
\fc  but it may not be so at  time $t=1$ if the smallest element of $\Ac^{cg,s}(0)$ is not in $\Ac^{er}(0)$.
\nc
Therefore we need 
a more subtle construction, proceeding with common elements as much as possible. We {\it delay}
 the elements $\Dc^{cg,d}$ corresponding to servers which are  informed in the {\it complete graph} dynamics but not yet informed on the Erd\"os-R\'enyi graph,  
 performing the construction with the common servers as much as possible. Thus the construction for the complete graph remains close to the one for the Erd\"os-R\'enyi graph.
\medskip

{\bf Delayed construction.} 
We construct dynamically the random tree $\Tc^{cg,d}(t)$ and the labelling $L=L^{cg,d}$
{\em simultaneously} with that on  the Erd\"os-R\'enyi graph according to (\ref{eq:er0}), (\ref{eq:ert}).
Initially, in addition to the sets defined in (\ref{eq:er0}), we also consider
\begin{align}
\Dc^{cg,d}(0)&=
\{w\in \Xc^{er}(0): L(w)\neq i, L(w)=L(w') \text{ for some } w'<w, w' \in \Xc(0)\setminus \Xc^{er}(0),\nonumber\\
&\qquad \qquad L(w)\neq L(w''), w''<w, w'' \in \Xc^{er}(0)\}
, \nonumber\\
\Ic^{cg,d}(0)&=\{\Par\}\cup \Ac^{er}(0)\cup \Dc^{cg,d}(0).
\label{eq:cgd0}
\end{align}
For the delayed dynamics on the complete graph, the partition in (\ref{eq:partition}) has three terms: $\Ec(0)=\{\Par\}, \Ac=\Ac^{er}(0)$
and $\Dc=\Dc^{cg,d}(0)$. The label function is defined in (\ref{eq:er0}).

With the process $(\Ic^{er}(t), \Ic^{cg,d}(t), \Ac^{er}(t),  \Dc^{cg,d}(t))$ and $L$ at time $t$, the value at the next step $t+1$ is defined by:

\begin{itemize}

\item If $\Ac^{er}(t)$ is non empty, we follow all the prescriptions in (\ref{eq:ert}), and we also define 
$X^{cg,d}(t)=X^{er}(t)$ denoted by $X(t)$ therein, 
\begin{align*}
C(t+1)&=\{w\in \Xc^{er}(t\!+\!1) :  L(w)\notin L(\Ic^{cg,d}(t)), L(w)= L(w')\text{ for some }w'\!<\!w, \\
& \qquad  w'\in \Xc(t\!+\!1)\setminus \Xc^{er}(t\!+\!1), L(w)\neq L(w''), w''<w, w'' \in \Xc^{er}(t+1)\}.
\end{align*}
(These are the nodes with label informed for the first time during the current burst on the complete graph, still not informed on the Erd\"os-R\'enyi graph. They will be 
placed in the set $\Dc^{cg,d}$ of delayed servers.) We also define the set 
\[F(t+1)=\{ w\in \Wc:Ê w\in \Dc^{cg,d}(t), L(w)\in L(\Ac^{er}(t+1))\}.\]
We then update the sets
\begin{align}
\Dc^{cg,d}(t+1)&=\{w\in \Wc: w\in \Dc^{cg,d}(t), L(w)\notin L(\Ac^{er}(t+1))\}\cup C(t+1),\nonumber\\
\Ic^{cg,d}(t+1)&=(\Ic^{cg,d}(t)\setminus F(t+1)) \cup \Ac^{er}(t+1)\cup \Dc^{cg,d}(t+1). \label{eq:2k10}
\end{align}
During this step, the mapping $L^{cg,d}=L^{er}=L$ is extended according to the rule in (\ref{eq:ert}).
The partitions (\ref{eq:partition}) for $\Tc^{er}(t)$ and $\Tc^{cg,d}(t)$ are then given by $\Ac=\Ac^{er}(t), \Ec=\Ec(t)$ (defined by $\Ec(t+1)=\Ec(t)\cup\{X(t+1)\}$) for both cases, and $\Dc=\emptyset$ for the first one but 
$\Dc=\Dc^{cg,d}(t)$ for the second one. It is therefore natural to set $\Ac^{cg,d}(t)=\Ac^{er}(t)$ for all such $t$'s. 
This step is in force till the first time when $\Ac^{er}(t)= \emptyset$, i.e. at time 
$\tau_n^{er}=t$ when the information process on the Erd\"os-R\'enyi graph ceases to evolve, after which  we proceed as follows. 
\item If $\Ac^{er}(t)$ is empty and $\Dc^{cg,d}(t)$ is non empty,  we set $X^{cg,d}(t+1)=X(t+1)$
with
\[
X(t+1)=\inf\{w\in \Dc^{cg,d}(t)\},\phantom{*}\mbox{denoted by $v$}
\]
and we obtain a realization of $K_{Lv}$, $\Ber^{Lv}_1,\dots, \Ber^{Lv}_{K_{Lv}}$ and $I^{Lv}_1,\dots,I^{Lv}_{K_{Lv}}$. To each descendant $vl$, $1\leq l \leq K_{Lv}$ we associate the label $I^{Lv}_{l}$. Then, we define
\begin{align*}
\Xc^{cg,d}(t+1)&=\{vl \in \Wc : 1\leq l \leq K_{Lv},\Ber^{Lv}_l=1\},\\
L(w) &= I^{Lv}_l, \qquad w=vl \in \Xc^{cg,d}(t+1) \qquad ({\rm with} \; L=L^{cg,d}).\nonumber
\end{align*}
We then update the sets 
\begin{align}
\Dc^{cg,d}(t+1)&=\Big(\Dc^{cg,d}(t)\setminus \{v\} \Big) \cup C(t+1),\nonumber\\
C(t+1) &=
\{w\in \Xc^{cg,d}(t \! +\!1) :  L(w)\!\notin \!L(\Ic^{cg,d}(t)), L(w)\!\neq \!L(w'), w'\!<\!w, w'\in \Xc^{cg,d}(t\!+\!1)\},\nonumber\\
\Ic^{cg,d}(t+1)&=\Ic^{cg,d}(t) \cup \Dc^{cg,d}(t+1). \label{eq:cgdt}
\end{align}
Recalling that $\Ac^{cg,d}(t)=\Ac^{cg,d}(\tau_n^{er})=\emptyset$ for $t > \tau_n^{er}$, we see that the complementary set
$$\Ec^{cg,d}(t+1)= \Ic^{cg,d}(t+1) \setminus \Dc^{cg,d}(t+1)$$ 
evolves like
$\Ec^{cg,d}(t+1)=\Ec^{cg,d}(t) \cup \{X(t+1)\}$.
 
\item If $\Ac^{er}(t)$ and $\Dc^{cg,d}(t)$ are empty,  the evolution is stopped and we denote by 
$\tau_n^{cg,d}$ the smallest such time $t$. \qed
 \end{itemize}

\begin{prop} \label{prop:cgd} We have the equality in law of the processes
\begin{align}
\label{eq:eqlawcg1}
\big( {\rm card} \, \Tc^{cg,s}(t),  {\rm card} \; \Ac^{cg,s}(t),  {\rm card} \, \Ec^{cg,s}(t) \big)_{t\geq 0} \qquad \qquad \qquad \qquad \qquad \qquad \qquad \qquad \qquad  \\ \eqlaw 
\big( {\rm card} \, \Tc^{cg,d}(t),  {\rm card} (\Ac^{cg,d}(t) \cup \Dc^{cg,d}(t)),  {\rm card} \, \Ec^{cg,d}(t) \big)_{t\geq 0}\;.
\quad  \nn
\end{align}
In particular,
\begin{equation}
\label{eq:eqlawcg2}
\big( {\rm card} \, \Tc^{cg,s}(\8),  \tau_n^{cg,s}\big) \eqlaw 
\big( {\rm card} \, \Tc^{cg,d}(\8),  \tau_n^{cg,d} \big)\;.
\quad  \nn
\end{equation}
\end{prop}
$\Box$ Since for all $t$, 
\begin{eqnarray*}
{\rm card} \, \Tc^{cg,s}(t)&=& {\rm card} \; \Ac^{cg,s}(t) + {\rm card} \, \Ec^{cg,s}(t),\\
 {\rm card} \, \Tc^{cg,d}(t)&=& {\rm card} (\Ac^{cg,d}(t) \cup \Dc^{cg,d}(t))+ {\rm card} \, \Ec^{cg,d}(t),
 \end{eqnarray*}
and for $i=s, d$, 
$$
 {\rm card} \, \Ec^{cg,i}(t)=(t \wedge \tau_n^{cg,i})+1, 
 $$
 it is enough to show that 
 \begin{equation}
\nn
 \big({\rm card} \; \Ac^{cg,s}(t); t \geq 0 \big) \eqlaw \big( {\rm card} (\Ac^{cg,d}(t) \cup \Dc^{cg,d}(t)); t \geq 0 \big).
\end{equation}
This relation follows from the equalities of the transitions
$$
\IP\Big(  {\rm card} \Ac^{cg,s}(t+1)= \cdot \mid {\rm card} \Ac^{cg,s}(s)=a_s, {\rm card} \Ec^{cg,d}(s)=s+1, s \leq t \big)
$$
and
$$
\IP\Big(  {\rm card}(\Ac^{cg,d}(t+1) \cup \Dc^{cg,d}(t+1))= \cdot \mid
 {\rm card}(\Ac^{cg,d}(s) \cup \Dc^{cg,d}(s))=a_s, {\rm card}\Ec^{cg,d}(s)=s+1, s \leq t \big)\;,
$$
for all $t \geq 0$. Indeed, from (\ref{eq:cgst}), and from (\ref{eq:ert}, \ref{eq:cgdt}),  both transitions
are equal to the law of the variable $a_t-1+Y$, with $Y$ the number of new coupons obtained 
in $\widehat K$ attempts in a coupon collector process with $n$ different coupons starting with
initially $a_t+t+1$
already obtained coupons.
 \qed
\medskip

Hence the sequential and delayed  constructions are equivalent to the standard transmission processs on the complete graph. Here is a direct consequence of  Propositions \ref{prop:cgs} and \ref{prop:cgd}.

\begin{cor}\label{cor:taucgsd}
 It holds $\tau_n^{cg,s} \eqlaw \tau_n^{cg,d}$ and $\widehat R( \tau_n^{cg,d}) \eqlaw {\mathfrak  T}_n$.
\end{cor}

In fact, we have a stronger result.

\begin{prop} \label{prop:2cg}
It holds  
$L(\Ic^{cg,s}(\infty))=L(\Ic^{cg,d}(\infty))$ for all $\omega$.
\end{prop}

$\Box$ This set depends only on the arrows, not on the order. 
Mathematically, for $i, j \in {\cal N}$, write $i \leadsto j$ if there exists $1 \leq k \leq K_i$ with
$B^i_k=1$ and $I^i_k=j$. Then, it is not difficult to see that
$L(\Ic^{cg,s}(\infty)) =L(\Ic^{cg,d}(\infty))$ because both are equal to the union
$$
\cup_{m=0}^{\8} \{i \in {\cal N}: \exists  i_0,\ldots i_m   \in {\cal N},
i_0=I_0, i_m=i, i_{k} \leadsto i_{k+1}, 0\!\leq \!k \!\leq \!m\!-\!1\},
$$
the above set being understood as  $ \{I_0\} $ for $m=0$.
\qed

\subsection{Coupling results}

With the above constructions, the information processes on the complete graph, in its delayed version, and on the 
Erd\"os-R\'enyi graph   are finely coupled.
First, it directly follows from the construction  that:
\begin{equation} \label{eq:couplage1}
\tau_n^{cg,d}\geq \tau_n^{er}, \qquad \Ic^{er}(t)\subset \Ic^{cg,d}(t)\quad {\rm for\ all\ } t\geq 0,   
\end{equation}
and 
\begin{equation} \label{eq:couplage1}
\Ac^{cg,d}(t)=\Ac^{er}(t)  \quad {\rm for\ all\ }   t\leq \tau^{er}_n.
\end{equation}

Note also that, from (\ref{eq:cgdt}),  the set $\Dc^{cg,d}(t)$ can increase or decrease with time, and that the elements of
$\Dc^{cg,d}(\tau_n^{er})$ together with their descendants encode the difference between the two 
processes.

\begin{prop} \label{prop:couplage-tau}
Assume $\IE  K^2 <\8$.
Then, we have 
\begin{equation}
\tau_n^{cg,d}-\tau_n^{er}=O_P(1),
\label{eq:taucgd=er}
\end{equation}
and
\begin{equation}
\label{eq:Tccgd=er}
{\tt card} \big( \Tc^{cg,d}(\8) \setminus  \Tc^{er}(\8) \big) = O_P(1).
\end{equation}
\end{prop}
We recall that for a sequence of real random variables $Z_n$, we write $Z_n=O_P(1)$  when 
$\sup_{n} \IP( |Z_n| \geq A) \to 0$ as $A \to +\8$. By Markov inequality,   a sufficient condition for that is $\sup_n \IE |Z_n| < \8$.
\medskip

$\Box$ First, observe that 
\begin{eqnarray}
\label{eq:j14-1}
\Tc^{cg,d}(t)&=& \Tc^{er}(t)\cup \Dc^{cg,d}(t), \qquad t \leq \tau_n^{er},\\
\label{eq:j14-2}
\Tc^{cg,d}(t)&=& \Tc^{er}(t)\cup \Dc^{cg,d}(t) \cup \big( \Ec^{cg,d}(t) \setminus \Ec^{cg,d}(\tau^{er}_n)\big),
 \qquad t \geq \tau^{er}_n.
\end{eqnarray}
For $t \leq \tau^{er}_n$,  in view of (\ref{eq:2k10}), the set $\Dc^{cg,d}(t)$ can increase at most by $C(t)$. Letting $i=LX(t)$, we observe that $\Dc^{cg,d}(t)$ is added a node with label $j \in \Nc$ if $j$
appears at least twice in $(I^i_k; k \leq K_i)$, first with a Bernoulli $B^i_k=0$ and then at least once 
 with a Bernoulli $B^i_k=1$. Then, for $i,j \in \Nc$,  we define the event $M(i,j)$ and the random variable $M(i)$
 \begin{eqnarray*}
M(i,j)&=&\big\{\exists k_1<k_2\leq K_i: B^i_{k_1}=0, B^i_{k_2}=1, I^i_{k_1}=I^i_{k_2}=j\big\},\\
M(i)&=& \sum_{j \in \Nc} \1{M(i,j)},
\end{eqnarray*}
and we have, from the above observation, 
\begin{equation}
\nn
{\rm card} \; \Dc^{cg,d}(t) -  {\rm card} \; \Dc^{cg,d}(t-1) \leq M(LX(t)).
\end{equation}
Thus,
\begin{equation}
\label{eq:majcle}
{\rm card} \; \Dc^{cg,d}(t) \leq \sum_{i \in \Nc} M(i) \stackrel{\rm def}{=} Y, \qquad
t \leq \tau^{er}_n,
\end{equation}
since each label $i \in {\cal N}$ can be picked at most once. The positive variable $Y$ has mean
\begin{eqnarray*}
\IE Y &=& n^2 \IE [\IP( M(i,j)| K_i)] \\
&\leq& n^2 \IE   {K_i \choose 2} p(1-p) \frac{1}{n^2}\\
&\leq& \frac{p(1-p)}{2} \IE K^2.
\end{eqnarray*}
Since $K$ is square integrable, this is bounded, and then 
\begin{equation}
\nn
{\rm card} \; \Dc^{cg,d}(t)=
O_{P}(1), \qquad t \leq \tau^{er}_n, 
\end{equation}
and then
\begin{equation}
{\tt card} \big( \Tc^{cg,d}(\tau^{er}_n) \setminus  \Tc^{er}(\tau^{er}_n) \big) = O_P(1),
\end{equation}
by (\ref{eq:j14-1}). Since $\Ac^{er}(\tau^{er}_n)=\emptyset$, we also have 
$
{\tt card} (\Ac^{cg,d}(\tau^{er}_n) \cup \Dc^{cg,d}(\tau^{er}_n) )
= O_P(1),$
and then 
\begin{equation}
\label{eq:sam1}
{\tt card} \;\Ac^{cg,s}(\tau^{er}_n) 
= O_P(1).
\end{equation}
We claim that this implies 
\begin{equation}
\label{eq:sam2}
\widehat S_n^{cg,s}(\widehat R(\tau^{er}_n) )
= O_P(1).
\end{equation}
Indeed, the conditional law of ${\tt card} \;\Ac^{cg,s}(t)$ given  
${\tt card} \;\Ac^{cg,s}(t-1)=a_{t-1},  \widehat S_n^{cg,s}(\widehat R(t))=m$
is the law of the variable $a_{t-1}-1+Y$, with $Y$ the number of new coupons obtained in 
$m$ attempts in a coupon collector process with $n$ different images starting initially with 
$a_{t-1}+t$ already obtained different coupons. Hence, for (\ref{eq:sam1}) to hold, 
it is necessary that  (\ref{eq:sam2})  holds.

We now use a lemma, which deals with the complete graph case only.
 \begin{lm} \label{lem:keycg} Consider the process on the complete graph defined in (\ref{eq:CDSdyn}).\\
 (i) For $A,B>0$ define
 $$
 u(A,B)= \limsup_{n \to \8} \IP( \inf \{ S_n(t) ; t \in [A,  {\mathfrak T}_n-A]\} \leq B),
 $$
 with the convention $\inf \emptyset = +\8$. Then, for all finite $B$, $u(A,B) \to 0$ as $A \to \8$.\\
 (ii) In particular, 
 for any random sequence $\sigma_n$, 
 \begin{equation}
\nn
S_n(\sigma_n) = O_P(1)  \implies \min \{ \sigma_n, ( {\mathfrak T}_n- \sigma_n )^+ \}=O_P(1).
\end{equation}
 \end{lm}
 $\Box$ Proof of Lemma \ref{lem:keycg}:  If $q=0$, we have $  {\mathfrak T}_n=O_P(1)$ and the result is trivial. We focus on the case 
 $q\in (0,1)$.
  The infimum is finite if and only if $  {\mathfrak T}_n\geq 2A$, so we get
 \begin{eqnarray}
\nn
\limsup_{n \to \8}  \IP\big( \inf \{ S_n(t) ; t \in [A,  {\mathfrak T}_n-A]\} < \8 ,  {\mathfrak T}_n<n q/2\big)
&=& 
\lim_{n \to \8}  \IP( {\mathfrak T}_n\geq 2A,  {\mathfrak T}_n<n q/2)
\\&=& \label{eq:s5-1}
\IP(\tau^{GW} \in [2A, \8))
\end{eqnarray}
where $\tau^{GW} $ is the survival time of the Galton-Watson process with offspring distribution $K$.
Thus, the last term vanishes as $A \to \8$. 
 
 We now study the contribution of the event $\{  {\mathfrak T}_n\geq n q/2\}$. 
 From the computations in the proof of Theorem 2.2  in \cite{CDS12}, the process $S_n$ increases linearly on the survival set with slope $\IE K-1>0$ for times $t, t \to \8, t = o(n)$. Precisely, we can fix  $\eta>0$ and $\delta>0$ such that
 \begin{equation}
\label{eq:s5-2}
 \liminf_{n \to \8}  \IP(  S_n(t) \geq \eta t; t \in [A, n\delta] \mid {\mathfrak T}_n\geq n q/2) = 1- \eps(A),
 \quad {\rm with} \; \lim_{A \to \8} \eps(A) =0.
\end{equation}
Similarly, from the  law of large numbers at times close to $n q$, we see that  the process $S_n$, on the survival set 
decreases linearly at such times with slope
$e^{-q} \IE K-1<0$.  Precisely, we can choose  $\eta$ and $\delta$ such that we have also
 \begin{equation}
\label{eq:s5-3}
 \liminf_{n \to \8}  \IP(  S_n(t) \geq \eta ( {\mathfrak T}_n-t); t \in [ {\mathfrak T}_n-n\delta,  {\mathfrak T}_n-A] \mid
  {\mathfrak T}_n\geq n q/2) = 1- \eps(A),
\end{equation}
with some function $\eps$ such that $ \lim_{A \to \8} \eps(A) =0.$
Finally,  from large deviations, we have
\begin{equation}
\label{eq:s5-4}
 \lim_{n \to \8}  \IP(  S_n(t) \geq C n  ; t \in [n\delta, n(q-\delta/2)] \mid
  {\mathfrak T}_n\geq n q/2) = 1.
\end{equation}
From (\ref{eq:s5-2}), (\ref{eq:s5-3}) and (\ref{eq:s5-4}), we conclude that
$$
\liminf_{n \to \8}  \IP(  S_n(t) \geq \eta A; t \in [A,  {\mathfrak T}_n-A] \mid  {\mathfrak T}_n\geq n q/2) = 1- \eps(A),
$$
with $\eps$ as above. This, in addition to (\ref{eq:s5-1}), implies our claim (i).  The other claim (ii) follows directly from (i). \qed
\medskip

With the lemma we complete the proof of Proposition \ref{prop:couplage-tau}.
From (\ref{eq:sam2}), the lemma  shows that 
$\widehat R(\tau^{er}_n) $ is close to 0 or to $\widehat R( \tn^{cg,d})$.
In turn this  implies, by definition of $\widehat R$,  that 
\begin{equation}
\nn
 \min \{ \tau^{er}_n, |\tau^{er}_n - \tau^{cg,d}_n|\}=O_P(1).
\end{equation}
Moreover, it is not difficult to see directly from the construction that for $\eps>0$ small enough 
(in fact, $\eps < \widehat q$), 
$$
\lim_{n \to \8} \IP( \tau^{er}_n \geq n \eps) = \lim_{n \to \8} \IP( \tau^{cg,d}_n \geq n \eps) =1- \widehat \sigma^{GW}.
$$

Together with $\tau^{er}_n \leq \tau^{cg,d}_n $, the last two relations  imply that $|\tau^{er}_n - \tau^{cg,d}_n| =O_P(1)$, which is
(\ref{eq:taucgd=er}).
\medskip

Further, following  
\cite{CDS12}, we see that the subtree generated by  $\Dc^{cg,d}(\tau_n^{er})$ is subcritical. 
Indeed, similarly to above (\ref{eq:s5-3}), from the  law of large numbers at times $s \sim n \widehat q$, we see that  the process $\widehat S_n^{cg,d}$, on the survival set 
decreases linearly at such times with slope
$e^{-\widehat q} p \IE K-1<0$.
By (\ref{eq:j14-2}), this yields the desired conclusion (\ref{eq:Tccgd=er}).

 \qed

From these estimates we derive our main results for the first mode of emission.
\medskip 

$\Box$ {\em Proof of  Theorems \ref{th:LLN} and \ref{th:TLC}. }
The estimates (\ref{eq:taucgd=er}) and (\ref{eq:Tccgd=er}) are good enough to
apply Theorem A for the complete graph and resource $\widehat K$. Indeed, by (\ref{eq:taucgd=er}) we have
$$
\tau_n^{er} = \tau_n^{cg,d} + O_P(1),
$$
and the sequence $(\tau_n^{cg,d})_{n \geq 1}$  obeys the law of large numbers 
in (i) and the central limit theorem in (ii) of  Theorem A. Then,  
$(\tau_n^{er})_{n \geq 1}$  obeys the same limit theorems.
\qed

\section{Second mode of emission}\label{sec:2frigos}
 In this second part we consider, a slightly different kind of emission on the Erd\"os -R\'enyi random graph. At time $0$ a server $i$ is chosen uniformly among the $n$ servers. Then, this server chooses uniformly a target server $i_1$ among the $n$ servers. If the edge between $i_1$ and $i$ is present then $i$ transmits the information to $i_1$ and wastes one unit of its resource. If the edge between $i$ and~$i_1$ is absent, nothing happens. This operation is repeated until $i$ exhausts all of its resource $K_i$. Then, we chose another informed server and repeat the same procedure as for $i$. The process ends when all the informed servers have exhausted their resources. This mode of emission differs from the previous one in the fact that a server can only use a unit of resource if the edge between it and its target server is present. Hence it is a 
 perturbation of the information transmission process on the complete graph with resource $K$, but not $\widehat K$ in contrast with the above case.

On the probability space $(\Omega, \mathcal{F}, \IP)$, we consider, as in beginning of  Section \ref{sec:constr}, the random elements
$(K_i)_{i\in \Nc}$, $I_0$, $(I^i_k)_{(i,k)\in \Nc\times \N}$ and  $(\Ber^i_{k})_{i \in {\cal N}, k \geq 1}$. 

\medskip

We first explain the ideas of the construction.
We attach to each edge $e \in {\cal E}$ a variable $B(e) \in\{0,1, {\rm ``unknown"}\}$ indicating the current status of the edge, i.e. if the edge is, respectively,  closed, open or still unknown at this stage.
The variables are updated during the construction, they start from $B(e)={\rm ``unknown"}$
and can turn from ${\rm "unknown"}$ to 0 or 1 when they appear.
Transmission on the complete graph occurs whenever meeting a Bernoulli variable $B^i_k$  
with the value 1, more precisely, at the $K_i$ smallest such $k$'s. 
In the  Erd\"os -R\'enyi  case, we check if the variable $B^i_k$  is compatible with the current status of the edge. Emissions with incompatibilities are placed in the set 
$\cal D$ of delayed elements, and will be recast afresh later in the case of the random graph. Compatible emissions are common to the two processes, and in the case of a previously uninformed target 
it is
placed in the set $\cal A$ and will be used in priority to emit in turn. When this set becomes empty, we end with the two sets ${\cal D}^{er}, {\cal D}^{cg}$,
that we process independently. All processes  being  delayed in the construction, we don't indicate it in the notations; Similarly,  we denote by $\Ac(t)$ the set of active vertices, since it is the same for the two processes. 
\medskip

Let $T_i=\inf\{j\geq 1: \sum_{l=1}^j\Ber^i_{l}=K_i\}$.
It is convenient  to initialize the process at time $t=-1$.
We start with $B(e)={\rm "unknown"} $ for all edge $e$,
with 
$$\Ac(-1)=\{\Par\}, \qquad 
\Ec^{er}(-1)=\Ec^{cg}(-1)=\Dc^{er}(-1)=\Dc^{cg}(-1)=\emptyset, \qquad 
L \Par= I_0.
$$

\medskip

With the process $(\Dc^{cg}(t), \Dc^{er}(t), \Ec^{cg}(t),\Ec^{er}(t), \Ac(t))$ at time $t$, its value at the next step $t+1$ is defined as follows (one can check that for $t=-1$ the first case below is in force
with $X(0)=\Par$): 
\begin{itemize}
\item If $\Ac(t)$ is non empty, we let $X(t+1)$ be its first element,
that we denote by $v$ for short notations, as well as $L(X(t+1))=i$.
We define the sets
\begin{align}
\Xc(t+1)&=\{vl \in \Wc : 1\leq l\leq T_i\},\nonumber\\
\Xc^{cg}(t+1)&=\{vl \in \Xc(t+1):  \Ber^i_l=1\}\nonumber , \\
L(vl) &= I^i_l, \qquad vl \in \Xc(t+1) \nonumber.
\end{align}

For edges $e$ of the form $e=\langle i, L(vk)  \rangle$ with $B(e)$ previously unknown, 
the value is being discovered, so we assign 
$$
B(e) = B^i_{\ell(e)} \qquad {\rm with} \; \ell(e)= \min\{ \ell: I^i_\ell=I^i_k\}.
$$ 
Define also 
\begin{align}
{\rm Inc}_\ell(t+1) &=\big\{ w \in \Xc (t+1): B^i_w \neq B(\langle i, Lw \rangle), B(\langle i, Lw \rangle) =\ell \big\}, \qquad \ell=0,1,  \nonumber\\
{\rm Inc}(t+1) &= {\rm Inc}_0(t+1) \cup {\rm Inc}_1(t+1).
 \nonumber
 \end{align}
 Here, ${\rm Inc}(t)$ is the set of incompatibilities at time $t$, they are of two possible nature. 
 With $[x]^+=\max\{x,0\}$, 
 \begin{equation}
\label{eq:canape2}
 m(t+1)=[ {\rm card(Inc}_0(t+1))-{\rm card(Inc}_1(t+1))]^+
 \end{equation}
is the number of emissions from server $i$  on the random graph 
 to be recast later. 
An incompatibility from the set ${\rm Inc}_0(t+1)$ corresponds to an emission 
 on the complete graph, but not on the Erd\"os-R\'enyi graph, and it is delayed. 
 An incompatibility from the set ${\rm Inc}_1(t+1)$ corresponds to an emission 
 on the Erd\"os-R\'enyi graph (but not on the complete graph). Note that this does not increase the number of informed servers.
 Define
\begin{equation}
\label{eq:canape}
T'_i= \max \left\{ \ell \leq T_i: \sum_{k=1}^\ell {\bf 1}_{vk \in {\rm Inc}_1(t+1)}
+  \sum_{k=1}^\ell {\bf 1}_{vk \in \Xc^{cg}(t+1) \setminus {\rm Inc}_0(t+1)} \leq K_i\right\}.
\end{equation}
 Then, we update 
 \begin{align}
\Ec^{cg}(t+1)&=\Ec^{er}(t+1)= \Ec^{er}(t) \cup \{v\},\nonumber\\
{\cal D}^{cg}(t+1)&= {\cal D}^{cg}(t) \cup {\rm Inc}_0(t+1)
\cup \big \{vk \in \Xc^{cg}(t  \!  + \! 1) \setminus   {\rm Inc}_0(t \! + \! 1) : k > T_i'\big\}, \nonumber\\
\Ac(t+1)&= (\Ac(t)\setminus \{v\}) \nonumber \\
&\!\!\!\! \cup
\Big \{vk \in \Xc^{cg}(t  \!  + \! 1) \setminus   {\rm Inc}_0(t \! + \! 1) : k \leq T_i', Lw \neq Lw',  w<w' \; {\rm and\ }  w' \in 
\Xc^{cg}(t \! + \! 1),
\nonumber\\
& \qquad \qquad \qquad
 Lw \neq Lw'',   w'' \in \Ec(t+1) \cup \Ac(t) \Big\},
\nonumber\\
{\cal D}^{er}(t+1)&={\cal D}^{er}(t) \cup  
\left\{ vk \in {\rm Inc}(t+1): \sum_{\ell=1}^k {\bf 1}_{v\ell \in {\rm Inc}(t+1)} \leq m(t+1)\right\}
 \label{eq:tr2-1}
\end{align}

\item When $\Ac(t)$ becomes empty, we set  $\widetilde \tau_n= t$, and from that time on, 
we continue  {\it separately}
the transmission processes on each of the two graphs, with the delayed emissions from the sets
$\Dc^{er}(\widetilde \tau_n), \Dc^{cg}(\widetilde \tau_n)$. They will terminate at later times $ \bar \tau_n^{er}, \tau_n^{cg}$ when $\Dc$ gets empty. 

(i) For the step from times $t$ to $t+1$ on the complete graph, we let $v$ be the first element  of $ \Dc^{cg}(t)$ and $i=Lv$. If the label $i$ is not an element of $L \Ec^{cg}(t)$, we update 
$\Ec^{cg}(t+1)=\Ec^{cg}(t) \cup \{v\}$  and $\Dc^{cg}(t+1)=(\Dc^{cg}(t)\setminus \{v\} ) \cup \{vk; k=1,\ldots K_i\}$. If
the label $i$ is an element of $L \Ec^{cg}(t)$, we update $\Ec^{cg}(t+1)=\Ec^{cg}(t)$  and $\Dc^{cg}(t+1)=\Dc^{cg}(t)\setminus \{v\} $. We then go to the next step.

When $\Dc^{cg}(t)$ becomes empty, set $\tau_n^{cg}=t$.

(ii) For the  Erd\"os-R\'enyi graph, for the step from times $t \geq \widetilde \tau_n$ to $t+1$:
\begin{itemize}
\item If there is a $s \leq \widetilde \tau_n$ such that ${\rm Inc}(s) \bigcap \Dc^{er}(t)
\neq \emptyset  $, 
consider the smallest one, still denoted by $s$, $v=X(s)$, $m(s)$ from (\ref{eq:canape2})
and $i=Lv$. Scan the edges
$e=\langle i, I^i_k\rangle$ for $k=T_i+1,T_i+2,\ldots$  to find the $m(s)$ first ones which are in the graph, and let $k_1,\ldots, k_{m(s)}$ the corresponding indices; If an edge $e$ was still unknown, 
put $B(e)=B^i_k$ for these $k$'s. Then, update 
$\Ec^{er}(t+1)=\Ec^{er}(t)$  and $\Dc^{er}(t+1)=\big(\Dc^{er}(t)\setminus  {\rm Inc}(s)\big)
\cup \{vk_1,\ldots, vk_{m(s)}\}.
$
Then go the next step.

\item If there is no $s \leq \widetilde \tau_n$ with 
 ${\rm Inc}(s) \bigcap \Dc^{er}(t)
\neq \emptyset  $, 
consider the smallest element $v$ in $\Dc^{er}(t)$ if any,    and $i=Lv$. Scan the edges
$e=\langle i, I^i_k\rangle$ for $k\geq 1$  to find the $K_i$ first ones which are in the graph, and let $k_1,\ldots, k_{K_i}$ the corresponding indices; If an edge $e$ was still unknown, 
put $B(e)=B^i_k$ for these $k$'s. Then, update 
$\Ec^{er}(t+1)=\Ec^{er}(t) \cup \{v\}$  and $\Dc^{er}(t+1)=\big(\Dc^{er}(t)\setminus  \{v\} \big)
\cup \{vk_1,\ldots, vk_{K_i}\}.
$
Then go the next step.

\item When $\Dc^{er}(t)$ becomes empty, set $\bar \tau_n^{er}=t$.

\end{itemize}

\end{itemize}

The above construction is indeed a fine coupling of transmission processes on the  Erd\"os-R\'enyi and the complete graphs. 
\medskip

$\Box$ Proofs of Theorems \ref{th:LLN2} and \ref{th:TLC2}.
To get an incompatibility it is necessary to pick twice the same edge in the construction, and 
to meet an event of the type
$$
\widetilde M(i,j; k_1,k_2)=\big\{
\langle i, I^i_{k_1} \rangle = \langle j, I^j_{k_2} \rangle , 
B^i_{k_1} \neq B^j_{k_2} \big\}.
$$ 
More precisely, the following events are equal,
$$
\{ {\rm Inc}(t+1) \neq \emptyset \} = \bigcup_{j \in {\cal N}, k_1 \leq T_{LX(t+1)}, k_2 \leq T_j}
\widetilde M(LX(t+1),j; k_1,k_2).
$$
The sets $\Dc^{er}, \Dc^{cg}$ increase only because of incompatibilites.
Similar to (\ref{eq:majcle}),  we can estimate, for times $t \leq \widetilde  \tau_n$,  the size of both sets by
\begin{equation}
\label{eq:majcle1}
 {\rm card} \; \Dc^{er}(t) 
\leq \sum_{i, j \in \Nc} \sum_{k_1 \leq T_i, k_2 \leq T_j}
{\bf 1}_{\widetilde M(i,j; k_1, k_2)} \stackrel{\rm def}{=} \widetilde Y, \qquad
{\rm card} \; \Dc^{cg}(t)  \leq \widetilde Y, 
\end{equation}
since each server $i \in {\cal N}$ can emit at most one burst. An elementary computation shows that the expectation of
$ \widetilde Y$ is bounded in $n$ as soon as $\IE K^2 <\8$. Then we obtain
\begin{equation}
\nn
{\rm card} \; \Dc^{cg}(\widetilde \tau_n) =
O_{P}(1), \qquad
{\rm card} \; \Dc^{er}(\widetilde \tau_n)
 =
O_{P}(1).
\end{equation}
Following the line of proof of Proposition \ref{prop:couplage-tau}
and  using Lemma \ref{lem:keycg},
we derive from the first above estimate that 
$$
\tau_n^{cg} - \widetilde \tau_n = O_P(1).
$$
Now, let us see that $\bar{\tau}_n^{er} - \widetilde \tau_n = O_P(1)$:
we first prove that at time $\widetilde \tau_n$, with high probability, for all $i\notin \Ec(\widetilde \tau_n)$ the number of edges $e=\langle i,j  \rangle$ such that $B(e)\neq {\rm``unknown"}$ is bounded from above by $ \sqrt{n}\ln n$.

Indeed, denoting by $E_n^i$ the number of ``known" edges adjacent to $i\in \Nc$ { at time} $\widetilde \tau_n$ we have 
\begin{equation*}
E_n^i {\bf 1}_{\{ i\notin \Ec(\widetilde \tau_n)\}}=  \sum_{j \in {\cal N}\setminus \{i\}}{\bf 1}_{\{B(\langle i,j  \rangle)\neq {\rm``unknown"}\},
 i\notin \Ec(\widetilde \tau_n), j \in \Ec(\widetilde \tau_n)\}}
\end{equation*}
where -- here as well as below -- the values of the random variables $(B(\langle i,j  \rangle))_{j\in \Nc}$ are taken {\it at time} $\widetilde \tau_n$. We estimate the second moment
\begin{align}
\IE \big[(E_n^i)^2 {\bf 1}_{\{i\notin \Ec(\widetilde \tau_n)\}}\big]
&=   \sum_{j \in {\cal N}\setminus \{i\}} \IP\big(B(\langle i,j  \rangle)\neq {\rm``unknown"};
 i\notin \Ec(\widetilde \tau_n); j \in \Ec(\widetilde \tau_n)\big)\nonumber\\
&\phantom{**}+  \sum_{j\neq j' \in {\cal N}\setminus \{i\}}
 \IP\big(B(\langle i,j  \rangle), B(\langle i,j'  \rangle)\neq {\rm``unknown"};
 i\notin \Ec(\widetilde \tau_n); j, j' \in \Ec(\widetilde \tau_n)\big)\nonumber
\\&\leq   \frac{n-1}{n} \frac{\IE K}{p} +    \frac{(n-1)(n-2)}{n^2} \frac{\IE K^2}{p^2}, 
\label{RTO}
\end{align}
bounding the second line by
\begin{align}
\lefteqn{ 
\IE \sum_{j\neq j' \in {\cal N}\setminus \{i\}}
 \IP\big(B(\langle i,j  \rangle), B(\langle i,j'  \rangle)\neq {\rm``unknown"};
 i\notin \Ec(\widetilde \tau_n); j,  j' \in \Ec(\widetilde \tau_n)\mid K_j, K_{j'} \big)} \phantom{************}
 \nonumber
\\ \leq  {(n\!-\!1)(n\!-\!2)} &
\IE \left[
 \IP\big(B(\langle 1,2  \rangle), B(\langle 1,3  \rangle)\neq {\rm``unknown"} \mid
 1\notin \Ec(\widetilde \tau_n); 2, 3 \in \Ec(\widetilde \tau_n); K_2, K_{3} \big)\right] \phantom{************}
  \nonumber
\\ \leq  (n\!-\!1)(n\!-\!2)&
\IE \left[
 \IP\big(B(\langle 1,2 \rangle) \neq {\rm``unknown"} \mid
 1 \notin \Ec(\widetilde \tau_n); 2  \in \Ec(\widetilde \tau_n); K_2 \big)^2\right] \phantom{************}
 \nonumber \\
 \leq  (n\!-\!1)(n\!-\!2)  & \frac{\IE K^2}{n^2p^2}, \nonumber
\end{align}
and a similar, even simpler, bound for the second line of  (\ref{RTO}).
Combining the union bound and Markov inequality we deduce
\begin{align*}
\IP\Big(\sup_{i\notin \Ec(\widetilde \tau_n)}E_n^i>\sqrt{n}\ln n\Big)
&= \IP\Big(\bigcup_{i\in \Nc}\{E_n^i>\sqrt{n}\ln n,i\notin \Ec(\widetilde \tau_n)\}\Big)\nonumber\\
&\leq n \IP(E_n^i {\bf 1}_{\{i\notin \Ec(\widetilde \tau_n)\}}>\sqrt{n}\ln n) \nonumber\\
&\leq \frac{\IE K+\IE K^2}{p^2 \ln^2 n}
\end{align*}
by (\ref{RTO}), which goes to 0 as $n\to \infty$.
We immediately obtain that at time $\widetilde \tau_n$,  with high probability, for each element of $\Dc^{er}(\widetilde \tau_n)$ the number of ``known" incident edges is $o(n)$. We deduce that with high probability, as $n\to \infty$, the edges chosen to generate the subtrees of the elements of $\Dc^{er}(\widetilde \tau_n)$ are ``unknown". Since $\tau_n^{cg} - \widetilde \tau_n = O_P(1)$, the number of informed servers at time $\widetilde \tau_n$ is of order $qn$ on the survival set $\{\sup_n \widetilde \tau_n=\infty\}$. Now, from the last two comments and by stochastic domination by a Galton-Watson process with offspring mean $\frac{(1-q)n}{n-o(n)}\IE(K)$, which is asymptotically smaller than $1$,  we can conclude that the subtrees generated by the elements of $\Dc^{er}(\widetilde \tau_n)$ are subcritical. 
On the other hand, on the extinction set $\{\sup_n \widetilde \tau_n<\infty\}$, the probability that ${\rm card}(\Dc^{er}(\widetilde \tau_n))>0$ goes to 0 as $n\to \infty$.
Hence,  
\begin{equation}
\label{eq:fin}
\bar{\tau}_n^{er} - \widetilde \tau_n = O_P(1).
\end{equation}
 Therefore, we conclude  
that 
$\tau_n^{cg}-\bar{\tau}_n^{er}=O_P(1)$.

From that point the proof is completely similar to that of Theorems \ref{th:LLN} and \ref{th:TLC}.
We will not repeat the details.
\qed

\section*{Acknowledgements}
The 
authors thank the 
French-Brazilian 
program \textit{Chaires Fran\c{c}aises dans l'\'Etat
de S\~ao Paulo} which supported the visit of F.C.\ to Brazil.
C.G.\ thanks FAPESP (2013/10101-9) for financial support.
S.P.\ and M.V. were partially supported by
CNPq (grants 300886/2008--0 and 301455/2009--0). 
The last three authors thank FAPESP (2009/52379--8)  
for financial support. F.C.~is partially supported by CNRS, UMR 7599 LPMA.
The authors thank a referee for his careful reading and 
many suggestions that permitted to  improve the paper.

\small{

}

\end{document}